\documentclass[]{article}
 
 \usepackage[margin=1.25in]{geometry}
 \usepackage{booktabs}  
 \usepackage{amsmath}
 \usepackage{amsthm}
 \usepackage{textcomp}  
 \usepackage{xcolor}
 \usepackage{lipsum}    
 \usepackage{longtable}
 \usepackage{authblk}
 \usepackage{lipsum}
 \usepackage{amsfonts}
 \usepackage{graphicx}
 \usepackage{epstopdf}
 \usepackage{algorithmic}
 \ifpdf
 \DeclareGraphicsExtensions{.eps,.pdf,.png,.jpg}
 \else
 \DeclareGraphicsExtensions{.eps}
 \fi
 
 \usepackage{amsopn}
 
 \usepackage[colorlinks=true,allcolors=blue]{hyperref}
 \usepackage{bm} 
 \usepackage{xcolor,comment}
 \usepackage{enumitem}
 \usepackage[utf8]{inputenc}
 \usepackage[T1]{fontenc}
 
 \usepackage{tikz}
 \usetikzlibrary{matrix,shapes,decorations.pathreplacing,fit,backgrounds}
 
 \newcommand{\id}{\mathrm d}
 \newcommand{\vc}{\mathbf}

 \renewcommand{\tilde}{\widetilde}
 \newcommand{\pard}[2]{\frac{\partial #1}{\partial #2}}

 \newcommand{\phat}{\hat{p}}

 \DeclareMathAlphabet\mathbfcal{OMS}{cmsy}{b}{n}

 \newcommand{\R}{\mathbb{R}}

  \newcommand{\btheta}{\pmb\theta}
 
 \newtheorem{thm}{Theorem}[section]

 \DeclareMathOperator*{\argmin}{\text{argmin}}
 
 \title{Fisher information and shape-morphing modes for solving the Fokker--Planck equation in higher dimensions}
 \author{William Anderson}
 \author{Mohammad Farazmand\thanks{Corresponding author's email address: \href{mailto:farazmand@ncsu.edu}{farazmand@ncsu.edu}}}
 \affil{Department of Mathematics, North Carolina State University,\\ 2311 Stinson Drive, Raleigh, NC 27695-8205, USA}
 \date{}
 
\begin{document}
 	\maketitle
 	
 	\begin{abstract}
	The Fokker--Planck equation describes the evolution of the probability density associated with a stochastic differential equation. As the dimension of the system grows, solving
	this partial differential equation (PDE) using conventional numerical methods becomes computationally prohibitive. Here, we introduce a fast, scalable, and interpretable method 
	for solving the Fokker--Planck equation which is applicable in higher dimensions. This method approximates the solution as a linear combination of shape-morphing Gaussians with time-dependent means and covariances. These parameters evolve according to the method of reduced-order nonlinear solutions (RONS) which ensures that the approximate solution stays close to the true solution of the PDE for all times. As such, the proposed method approximates the transient dynamics as well as the equilibrium density, when the latter exists.
	Our approximate solutions can be viewed as an evolution on a finite-dimensional statistical manifold embedded in the space of probability densities. We show that the metric tensor in RONS coincides with the Fisher information matrix on this manifold. We also discuss the interpretation of our method as a shallow neural network with Gaussian activation functions and time-varying parameters. In contrast to existing deep learning methods, our method is interpretable, requires no training, and automatically ensures that the approximate solution satisfies all properties of a probability density.
\end{abstract}
 
 \section{Introduction}
Unlike deterministic dynamical systems, the evolution of stochastic systems cannot be unambiguously described based solely on the initial condition.
Instead one studies the probability distribution of the state over time~\cite{Sobczyk1991}. In most problems of practical interest, the evolution of the probability distribution cannot be determined analytically and 
therefore it needs to be approximated numerically.

For stochastic differential equations (SDEs), one can run a large ensemble of numerical simulations with different initial conditions and different realizations of the noise.
Subsequently, the probability distribution of the system state can be estimated using these large scale Monte Carlo simulations~\cite{silverman1986}.
Alternatively, one can numerically solve the Fokker--Planck equation which is a partial differential equation (PDE) describing the evolution of the probability density associated with the state of the system~\cite{pichler2013}.

The computational cost associated with both these approaches becomes quickly prohibitive as the dimension $d$ of the system grows~\cite{Chen2018,Mendez2022}. For instance, direct Monte Carlo methods need $\mathcal O(\epsilon^{-d})$ samples to reach an error tolerance $0<\epsilon\ll 1$ ~\cite{Naaman2021,Nagler2016,Stone1980}. On the other hand, discretizing the Fokker--Planck equation requires $\mathcal O(N^d)$ collocation points, where $N$ is the number of points in each direction.
In either case, the computational cost grows exponentially with the dimension of the system. 

We note that, for SDEs with special properties, there exist tailored methods with manageable computational cost in higher dimensions. These include SDEs with slow-fast dynamics~\cite{Majda2006} and Hamiltonian SDEs in equilibrium~\cite{Soize1988}.
We refer to~\cite{Chen2018} for a review of these special cases. Nonetheless, these methods are not applicable to general SDEs and often only approximate the equilibrium density, not the transient dynamics.

The excessive computational cost of solving PDEs in higher dimensions is not specific to the Fokker--Planck equation; discretizing any PDE in higher dimensions is computationally prohibitive. Only recently, deep learning methods have been able to overcome this curse of dimensionality~\cite{Weinan2018,Spiliopoulos2018}. In this approach, the solution of the PDE is approximated by a deep neural network. The parameters of the network are trained so that its output solves the PDE approximately. Being mesh-free, deep learning methods are better suited for solving PDEs in higher dimensions.

There is a rapidly growing list of such deep learning methods. For instance, physics informed neural networks (PINNs) train a deep neural network by minimizing the residual of the error at prescribed collocation points~\cite{Raissi2019}. 
Deep Galerkin method (DGM) takes a similar approach but instead of using collocation points, it minimizes the functional norm of the error~\cite{Spiliopoulos2018}. Consequently, since DGM does not use collocation points, it is specially suitable for solving PDEs in higher dimensions. Another notable example is neural operators~\cite{anandkumar2023} which learn maps between infinite-dimensional Banach spaces and can be used to solve PDEs~\cite{Anandkumar2022}.  We refer to Beck et al.~\cite{Beck2023} for a recent review of deep learning methods for solving PDEs.
Deep neural networks have already been used to solve the Fokker--Planck equation~\cite{beck2021,Chen2021,Tang2022,Kurths2020}. In spite of their impressive capabilities, these deep learning methods suffer from limited interpretability; the neural network is a black box, mapping initial-boundary conditions to the PDE's solution. 

Here, we introduce an alternative method based on reduced-order nonlinear solutions (RONS). RONS approximates the solution of a PDE as a linear combination of shape-morphing modes~\cite{anderson2021,anderson22,anderson2023}. In contrast to existing spectral methods, where the modes are static in time, RONS allows the modes to change shape and hence adapt to the solution of the PDE. This is achieved by allowing the modes to depend nonlinearly on a set of time-dependent shape parameters. The optimal evolution of the parameters are determined by solving a system of ordinary differential equations (ODEs), known as the RONS equation. In the case of the Fokker--Planck equation,  we choose Gaussians as our shape-morphing modes. The corresponding shape parameters are the mean and covariance of each mode which are allowed to change over time in order to better approximate the solution of the PDE.

As we discuss in Section~\ref{sec:modes} below, our method can be interpreted as a shallow neural network. However, it has several advantages compared to existing deep learning methods:
\begin{enumerate}
	\item Interpretability: Our method uses a linear combination of Gaussians with time-varying means and covariances. As such, the approximate solution can be easily interpreted in terms of the probability distribution of the system state. Further, a posteriori analysis of the solution is straightforward.
	\item No training required: The parameters of our solution evolve according to known and computable ODEs. As a result, no training (i.e., numerical optimization) and no data are needed for determining the parameters of the network.
	\item Conservation of probability: In RONS, it is straightforward to ensure that the approximate solution respects the conserved quantities of the PDE. In the case of the Fokker--Planck equation, this conserved quantity is the total probability, i.e., the integral of the probability density over the entire state space. As a result, our solutions are guaranteed to satisfy the properties of a probability density function.
\end{enumerate}
Our method can easily be extended to be used with deep neural networks~\cite{du21}. However, we intentionally use its shallow version to maintain its interpretability and low computational cost.
The universal approximation theorem of Park and Sandberg~\cite{Park1991} guarantees that probability densities can be approximated with such shallow neural networks to any desirable accuracy.

The remainder of this paper is organized as follows. In Section~\ref{sec:prelim}, we review the necessary mathematical preliminaries, describe the shape-morphing solutions which approximate the Fokker--Planck equation, and discuss their interpretation as a shallow neural network. Section~\ref{sec:evolution} reviews RONS for the optimal evolution of the shape parameters. In Section~\ref{sec:metric}, we discuss the relationship between our method and the Fisher information metric. The performance of our method is demonstrated on several numerical examples in Section~\ref{sec:results}. Finally, we present our concluding remarks in Section~\ref{sec:conc_FP}.

\section{Mathematical Preliminaries}\label{sec:prelim}
Consider the It\^o stochastic differential equation,
\begin{equation}
\id \vc X = \vc F(\vc X,t)\id t + \sigma\id \vc W,\quad \vc X(0)=\vc X_0
\label{eq:SDE}
\end{equation}
where $\vc X(t)\in\mathbb R^d$ denotes the state vector at time $t\geq 0$, $\vc F:\mathbb R^d\to\mathbb R^d$ is a sufficiently smooth vector field, and $\vc W(t)$ is the standard Wiener process in $\mathbb R^d$ with intensity $\sigma>0$. The initial condition $\vc X_0$ can itself be uncertain and drawn randomly from a probability density $p_0(\vc x)$.
Although here we only consider the homogeneous additive noise, the following framework can easily be extended to the case of inhomogeneous multiplicative noise where the noise intensity matrix $\sigma(\vc X,t)\in\mathbb R^{d\times d}$ depends on the state.

The probability density $p(\vc x,t)$ corresponding to the SDE~\eqref{eq:SDE} satisfies the Fokker--Planck equation,
\begin{equation}
\pard{p}{t} = \mathcal L p := - \nabla\cdot \left[ \vc F p\right]+\nu \Delta p,\quad p(\vc x,0)=p_0(\vc x),
\label{eq:FP}
\end{equation}
where the diffusion coefficient is given by $\nu = \sigma^2/2$ and $\mathcal L$ is a linear differential operator which depends on the vector field $\vc F(\vc x,t)$.
Being a probability density, the solution $p$ is non-negative and belongs to the Lebesgue space $L^1(\mathbb R^d)$~\cite{Otto1998}. Furthermore, the norm of the solution in this space is conserved and equal to unity, i.e.,
\begin{equation}
\| p(\cdot,t)\|_{L^1} = 1,
\end{equation}
for all time $t\geq 0$.

As discussed in the Introduction, when the dimension $d$ is large, discretizing the Fokker--Planck equation becomes computationally prohibitive. To overcome this curse of dimensionality, we use a mesh-free method by considering an approximate solution of the form,
\begin{equation}
\hat p(\vc x,\btheta(t)) = \sum_{i=1}^r A^2_i(t) \exp \left[ - \frac{|\vc x-\vc c_i(t)|^2}{L_i^2(t)}\right],
\label{eq:gen_approx}
\end{equation}
which is a sum of Gaussians. Here $A_i^2(t)\in\mathbb R^+$ is the amplitude of the $i$-th mode, $L_i^2(t)\in\mathbb R^+$ is its variance, and $\vc c_i(t)\in\mathbb R^d$ is its mean. The collection of amplitudes, variances and means constitutes the time-dependent shape parameters $\btheta = \left\{A_i,L_i,\vc c_i \right\}_{i=1}^r$.  Therefore, the approximate solution contains a total of $n=r(d+2)$ parameters.

The key to the success of our method is the fact that the parameters $\btheta(t)$ are allowed to change over time. This enables the modes in the approximate solution~\eqref{eq:gen_approx} to change their shape and location, hence adapting to the solution of the PDE. Of course, the immediate question is how to evolve the parameters $\btheta(t)$. As we review in Section~\ref{sec:evolution}, RONS 
evolves the parameters according to a set of ordinary differential equations. These ODEs are designed such that the solution $\hat p(\vc x,\btheta(t))$ approximates the dynamics of the Fokker--Planck equation. But before describing RONS, we first provide the justification for using Gaussian modes in the approximate solution~\eqref{eq:gen_approx}.

\subsection{Choice of the modes}\label{sec:modes}
The approximate solution~\eqref{eq:gen_approx} consists of a sum of Gaussians. In general, other nonlinear functions can be used as modes~\cite{anderson2021}. However, the Gaussian seems appropriate for the Fokker--Planck equation. First, note that if the vector field $\vc F(\vc x,t)$ is linear in $\vc x$, the stationary solution to the corresponding Fokker--Planck equation will be a Gaussian~\cite{Sobczyk1991}. More importantly, the following universal approximation theorem guarantees that any function in $L^1(\mathbb R^d)$ can be approximated, to any desirable accuracy,  with a function of the form~\eqref{eq:gen_approx}.

\begin{thm}[Park and Sandberg~\cite{Park1991}]\label{thm:UAT}
	Let $K:\mathbb R^d\to\mathbb R$ be integrable, bounded and continuous almost everywhere. If $\int_{\mathbb R^d}K(\vc x)\id \vc x\neq 0$, then the set 
	\begin{equation}
	\left\{\sum_{i=1}^r A_i K \left( - \frac{|\vc x-\vc c_i|^2}{L_i^2}\right): A_i\in\mathbb R, L_i\neq 0, \vc c_i\in\mathbb R^d, r\in\mathbb N \right\}
	\end{equation}
	is dense in $L^q(\mathbb R^d)$ for all $1\leq q<\infty$.
\end{thm}

The Gaussian function clearly satisfies all the conditions in Theorem~\ref{thm:UAT}. Therefore, any probability density function $p$ in $L^1(\mathbb R^d)$ can be approximated, to arbitrary accuracy, with a function $\hat p$ of the form~\eqref{eq:gen_approx}. More precisely, given $p\in L^1(\mathbb R^d)$, for any $\epsilon>0$, there exist shape parameters $\btheta$ and $r\in\mathbb N$ such that $\|p - \hat p\|_{L^1}<\epsilon$. We point out that the amplitudes $A_i$ in the approximate solution~\eqref{eq:gen_approx} are squared to ensure that $\hat p$ is non-negative as is required for a probability density.

Note that the covariance matrix of each Gaussian in~\eqref{eq:gen_approx} is diagonal. One could alternatively choose a non-diagonal covariance matrix $\Sigma_i(t)\in\mathbb R^{d\times d}$ and use the modes,
\begin{equation}
\exp\left[ -\frac12 (\vc x-\vc c_i)^\top \Sigma_i^{-1}(t)(\vc x-\vc c_i)\right].
\end{equation}
In higher dimensions, this would significantly increase the number of shape parameters per mode since the covariance matrix needs to be solved for simultaneously. Fortunately, Theorem~\ref{thm:UAT} implies that diagonal covariance matrices are sufficient for approximating any probability density as long as the number of modes $r$ is large enough.

\begin{figure}
	\centering
	\includegraphics[width=0.95\textwidth]{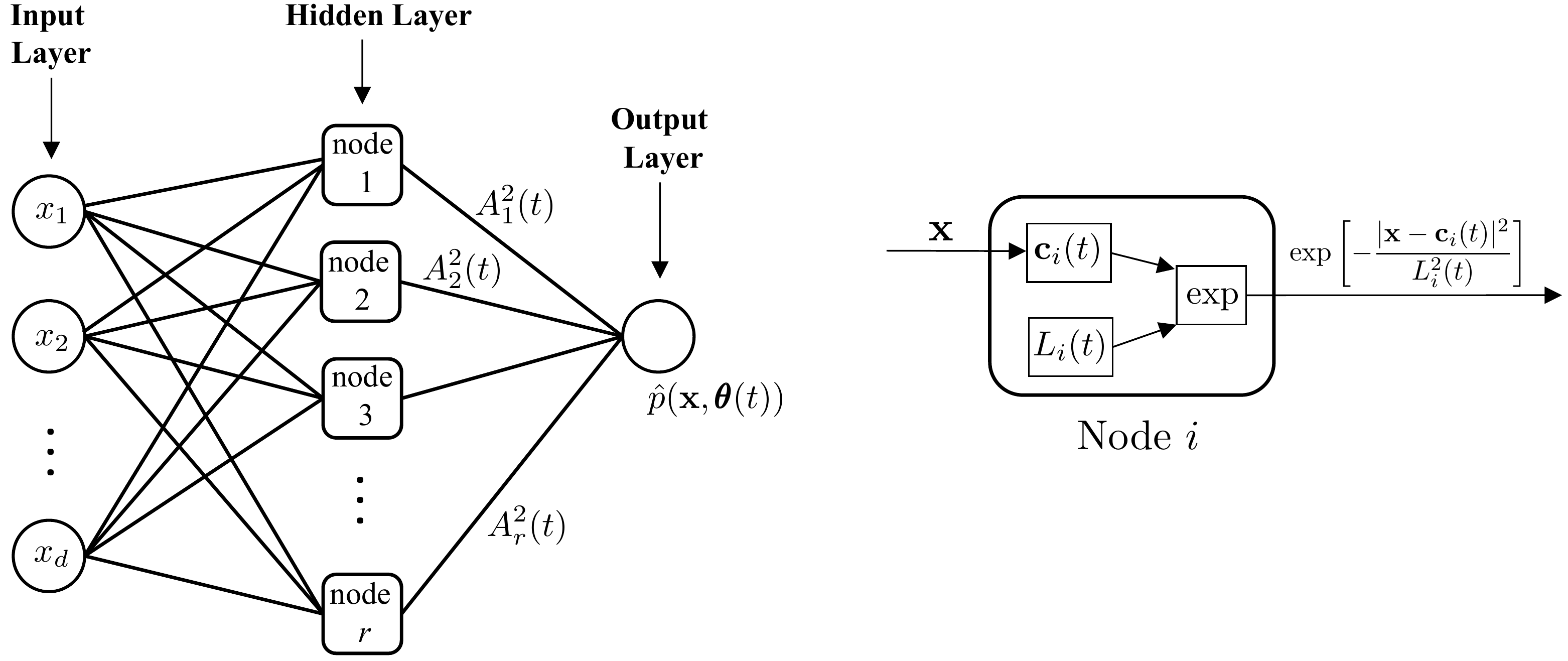}
	\caption{Interpretation of the approximate solution~\eqref{eq:gen_approx} as a shallow neural network. Left: Network architecture. Right: Internal structure of each node.}
	\label{fig:NN}
\end{figure}

We point out that the approximate solution~\eqref{eq:gen_approx} can be viewed as a shallow artificial neural network with one hidden layer. As depicted in figure~\ref{fig:NN}, this network takes the 
coordinates $\vc x = (x_1,x_2,\cdots,x_d)$ as inputs and returns $\hat p(\vc x,\btheta(t))$ as output. The activation function of each node is a Gaussian with parameters $L_i(t)$ and $\vc c_i(t)$. The amplitudes $A_i^2(t)$ are the output weights of the neural network. This network differs from conventional neural networks in that its parameters are time-dependent. As such, it is a special type of an evolutional neural network first introduced in~\cite{du21}.

Furthermore, in conventional neural networks, the parameters are determined through the so-called training process where the parameters are iteratively tuned to match the training data. 
In contrast, here we use RONS to evolve the parameters so that the approximate solution $\hat p(\vc x,\btheta(t))$ matches the dynamics of Fokker--Planck as closely as possible. This method, which requires no training and no data, is described in Section~\ref{sec:evolution} below.

\section{Evolution of Parameters}\label{sec:evolution}
RONS evolves the shape parameters $\btheta(t)$ such that the discrepancy between the approximate solution~\eqref{eq:gen_approx} and the dynamics of the PDE~\eqref{eq:FP} is instantaneously minimized~\cite{anderson2021}. Here, we briefly review this method in the context of Fokker--Planck equation and refer to Anderson and Farazmand~\cite{anderson2021,anderson2023} for further detail. 

Note that since $\hat p$ is an approximate solution $\partial_t\hat p$ does not necessarily coincide with the right-hand side $\mathcal L\hat p$ of the Fokker--Planck equation.
We define the residual
\begin{equation}
R(\vc x; \btheta,\dot{\btheta})=\sum_{j=1}^n \pard{\hat p}{\theta_j}(\vc x,\btheta)\dot{\theta}_j - \mathcal L\hat p(\vc x,\btheta),
\label{eq:res}
\end{equation}
where we used the fact that
\begin{equation}
\pard{}{t}\hat p(\vc x,\btheta(t))=\sum_{j=1}^n\pard{\hat p}{\theta_j}(\vc x,\btheta)\dot{\theta}_j.
\end{equation}

RONS determines the evolution of the parameters $\btheta(t)$ by minimizing the residual $R$ in an appropriate sense. We consider two options: 1) Symbolic RONS which minimizes a functional norm of the residual, and 2) Collocation RONS which minimizes the residual at prescribed collocation points. For completeness, we review each approach in Sections~\ref{sec:SRONS} and~\ref{sec:CRONS} below.

\subsection{Symbolic RONS}\label{sec:SRONS}
For a fixed $\btheta$ and $\dot{\btheta}$, consider $R(\vc x;\btheta,\dot{\btheta})$ as a map from $\vc x$ to the real line, i.e., $R(\cdot;\btheta,\dot{\btheta}):\mathbb R^d\to \mathbb R$.
We assume that this map belongs to a Hilbert space $H$ with the inner product $\langle \cdot,\cdot\rangle_H$ and the induced norm $\|\cdot\|_H$. Here, we describe RONS for a general Hilbert space and in Section~\ref{sec:metric} identify specific choices of the Hilbert space suitable for the Fokker--Planck equation.

In symbolic RONS, we minimize the residual norm $\|R(\cdot;\btheta,\dot{\btheta})\|_H$ with the constraint that 
\begin{equation}
I(\btheta) := \|\hat p(\cdot,\btheta)\|_{L^1} =1. 
\end{equation}
For the approximate solution~\eqref{eq:gen_approx}, we have $I(\btheta) = \sum_{i=1}^r A_i^2(\pi L_i^2)^{d/2}$.
This constraint is required to ensure that the approximate solution $\hat p$ is a probability density function. The resulting constrained optimization problem reads,
\begin{subequations}\label{eq:opt_srons}
	\begin{equation}\label{eq:min_srons}
	\min_{\dot{\btheta}} \|R(\cdot;\btheta,\dot{\btheta})\|_H^2+\alpha|\dot{\btheta}|^2,
	\end{equation}
	\begin{equation}
	\mbox{such that}\quad I(\btheta) = 1,
	\end{equation}
\end{subequations}
where $\alpha\geq 0$ is a Tikhonov regularization parameter and $|\cdot|$ denotes the usual Euclidean norm.
In the absence of regularization ($\alpha=0$), equation~\eqref{eq:min_srons} minimizes the instantaneous discrepancy between the rate of change of the approximate solution $\partial_t\hat p$ and the rate of change $\mathcal L\hat p$ dictated by the PDE. The motivation for the Tikhonov regularization will become clear in Section~\ref{sec:srons_comp}. One may also formulate a finite-time version of the constrained optimization problem~\eqref{eq:opt_srons}. However, this finite-time formulation tends to return unstable equations for the evolution of parameters $\btheta(t)$ (cf. Appendix A of \cite{anderson2021}).

As shown in~\cite{anderson2021,anderson2023}, the solution to the constrained optimization problem~\eqref{eq:opt_srons} satisfies the system of ordinary differential equations,
\begin{equation}\label{eq:srons}
\left[M(\btheta)+\alpha\mathbb I_n\right]\dot{\btheta} = \vc f(\btheta) - \lambda \nabla I(\btheta),
\end{equation}
where the gradient is taken with respect to the parameters $\btheta$, and $\mathbb I_n$ denotes the $n\times n$ identity matrix with $n=r(d+2)$ being the total number of parameters. 
We refer to equation~\eqref{eq:srons} as symbolic RONS equation, or S-RONS for short. The motivation for the term \emph{symbolic} will become clear at the end of this section.

The metric tensor $M(\btheta)$ is a symmetric positive semi-definite matrix whose entries are given by 
\begin{equation}\label{eq:srons_M}
M_{ij} = \left\langle \pard{\hat p}{\theta_i},\pard{\hat p}{\theta_j}\right\rangle_H,\quad i,j\in\{1,2,\cdots,n\}.
\end{equation}
Note that, since the metric tensor is symmetric positive semi-definite, the matrix $M(\btheta)+\alpha\mathbb I_n$ is invertible for all $\alpha>0$.
The vector field $\vc f:\mathbb R^n\to\mathbb R^n$ is defined by
\begin{equation}\label{eq:srons_f}
f_i = \left\langle \pard{\hat p}{\theta_i},\mathcal L\hat p\right\rangle_H,\quad i=1,2,\cdots,n.
\end{equation}
Finally, $\lambda\in\mathbb R$ is a Lagrange multiplier defined by
\begin{equation}
\lambda = \frac{\langle \nabla I,(M+\alpha\mathbb I_n)^{-1}\vc f\rangle}{\langle \nabla I,(M+\alpha\mathbb I_n)^{-1}\nabla I\rangle},
\end{equation}
where $\langle\cdot,\cdot\rangle$ denotes the usual Euclidean inner product. 

As shown in~\cite{anderson2023}, in the absence of regularization ($\alpha=0$), the ODEs~\eqref{eq:srons} become stiff as the number of parameters $n$ grows. As a result, solving the ODEs using explicit integration schemes requires exceedingly small time steps. However, as we show in Section~\ref{sec:results} below, even small values of the regularization parameter $\alpha>0$ alleviate this stiffness issue.

\subsubsection{Symbolic computing for S-RONS}\label{sec:srons_comp}
Finally, we comment on the computational cost of the symbolic RONS equation~\eqref{eq:srons}. This equation requires computing functional inner products $\langle\cdot,\cdot\rangle_H$ over the Hilbert space $H$. When feasible, we use symbolic computing to obtain closed-form symbolic expressions for the metric tensor $M_{ij}$ and the right-hand side vector field $f_i$; hence the name symbolic RONS.  This is desirable since symbolic expressions evade quadrature errors. Furthermore, the inner products do not need to be recomputed during time integration; the existing symbolic expressions are evaluated by substituting the updated values of the parameters $\btheta(t)$.

However, as the number of parameters $n=r(d+2)$ increases, brute-force computation of the inner products becomes expensive.
Computing the metric tensor, for instance, requires $\mathcal O(n^2)$ symbolic computations. Fortunately, as discussed in~\cite{anderson2023}, the special structure of the metric tensor 
implies that a far smaller number of symbolic computations are required. More specifically, the entire metric tensor $M\in\mathbb R^{n\times n}$ can be evaluated by computing only $6$ symbolic expressions. Similarly, the right-hand side vector field $\vc f\in\mathbb R^n$ can be evaluated by computing only $3$ symbolic expressions. 
Therefore, to form the entire S-RONS equation~\eqref{eq:srons}, only $9$ symbolic computation of inner products is needed. Note that this number is independent of the number of terms $r$ in the approximate solution~\eqref{eq:gen_approx} and the dimension $d$ of the SDE. As such, S-RONS is scalable to higher dimensions and the number of terms can be increased to achieved a desired accuracy at no significant additional computational cost.

To better understand this low computational cost, note that the inner products in the metric tensor~\eqref{eq:srons_M} involve the terms,
\begin{subequations}\label{eq:pard}
	\begin{equation}
	\pard{\phat}{A_k} = 2A_k\exp\left[-\frac{|\vc x-\vc c_k|^2}{L_k^2} \right],
	\end{equation}
	\begin{equation}
	\pard{\phat}{L_k} = \frac{2A_k^2}{L_k^3}|\vc x-\vc c_k|^2\exp\left[-\frac{|\vc x-\vc c_k|^2}{L_k^2} \right],
	\end{equation}
	\begin{equation}
	\pard{\phat}{c_{k,\ell}} = \frac{2A_k^2}{L_k^2}\left(x_\ell- c_{k,\ell}\right)\exp\left[-\frac{|\vc x-\vc c_k|^2}{L_k^2} \right],
	\end{equation}
\end{subequations}
where $1\leq k\leq r$, $1\leq \ell\leq d$, and $c_{k,\ell}$ denotes the $\ell$-th component of the vector $\vc c_k$. Assume that we have computed symbolic expressions for the inner products $\langle\cdot,\cdot\rangle_H$ between the terms in equation~\eqref{eq:pard} with general indices $(k,\ell)$. Then the entire metric tensor $M$ can be evaluated by substituting the numerical values of $A_k$, $L_k$, and $c_{k,\ell}$ in the computed symbolic expressions. Therefore, only 6 symbolic expressions need to be computed to populate the entire matrix $M$.

Similarly, to compute the right-hand side vector~\eqref{eq:srons_f}, we only need to compute symbolic expressions for the inner product of $\mathcal L\phat$ with the three terms in equation~\eqref{eq:pard}. Therefore, the entire vector $\vc f(\theta)$ can be evaluated using 3 symbolic expressions.

This observation leads to enormous computational savings as the number of dimensions and/or the number of terms in the approximate solution increase.
For instance, in $d=5$ dimensions and with $r=10$ terms in the approximate solution, brute-force computation of the S-RONS equation would require symbolic computation of 
$n(n+3)/2=2,555$ expressions, taking into account that the metric tensor $M$ is symmetric. 
However, as discussed above, in reality only $9$ symbolic expressions need to be computed in symbolic RONS. 

We emphasize that symbolic computing reduces the computational cost of the time integration as well. Since all terms are computed symbolically, as time integration progresses, these terms do not need to be recomputed; instead, they will be evaluated by substituting the new values of $\btheta(t)$ into existing symbolic expressions.

For certain choices of the Hilbert space $H$, existing symbolic computing software are unable to return a closed-form expression for the metric tensor or the right-hand side vector field. 
Collocation RONS addresses this issue.

\subsection{Collocation RONS}\label{sec:CRONS}
Symbolic computation of the RONS terms in~\eqref{eq:srons} may not be straightforward for certain Hilbert spaces $H$ and certain vector fields $\vc F(\vc x,t)$. 
Collocation RONS, or C-RONS for short, was developed in~\cite{anderson2023} to address this issue. In this approach, no functional inner products need to be computed and therefore the computational cost of forming the RONS equations reduces drastically. 

In C-RONS, instead of minimizing the residual $R$ over the entire state space $\mathbb R^d$, we minimize it over a set of prescribed collocation points $\{\vc x_1,\vc x_2,\cdots,\vc x_N\}$.
More precisely, we solve the constrained optimization problem,
\begin{subequations}\label{eq:opt_CRONS}
	\begin{equation}
	\min_{\dot{\btheta}} \sum_{i=1}^N|R(\vc x_i;\btheta,\dot{\btheta})|^2+\alpha|\dot{\btheta}|^2,
	\end{equation}
	\begin{equation}
	\mbox{such that}\quad I(\btheta) = 1.
	\end{equation}
\end{subequations}
In the absence of regularization ($\alpha=0$), this optimization problem minimizes the sum of squares of the residual at the collocation points given the constraint that $I(\btheta)=\|\hat p(\cdot,\btheta)\|_{L^1}=1$. As in S-RONS, the regularization ensures that the resulting ODEs are not stiff and therefore can be integrated in time using explicit discretization schemes.

As shown in~\cite{anderson2023}, the solution to the optimization problem~\eqref{eq:opt_CRONS} satisfies the system of ODEs,
\begin{equation}\label{eq:crons_ode}
\left[\tilde M(\btheta)+\alpha \mathbb I_n\right]\dot{\btheta} = \tilde{\vc f}(\btheta)-\tilde\lambda \nabla I(\btheta),
\end{equation}
where the collocation metric tensor is given by 
\begin{equation}
\tilde M(\btheta) = J^\top J, \quad J_{ij}(\btheta) = \pard{\hat p}{\theta_j}(\vc x_i,\btheta),\quad i\in\{1,2,\cdots,N\},\quad j\in\{1,2,\cdots,n\}.
\end{equation}
The right-hand vector field $\tilde{\vc f}:\mathbb R^n\to\mathbb R^n$ is given by 
\begin{equation}
\tilde{\vc f}(\btheta) = J^\top \vc f(\btheta),\quad f_i(\btheta)= \mathcal L\hat p\Big|_{\vc x=\vc x_i},\quad i\in\{1,2,\cdots,N\}
\end{equation}
and the Lagrange multiplier is defined by 
\begin{eqnarray}
\tilde\lambda = \frac{\langle \nabla I,[\tilde M(\btheta)+\alpha \mathbb I_n]^{-1}\tilde{\vc f}\rangle}{\langle \nabla I,[\tilde M(\btheta)+\alpha \mathbb I_n]^{-1}\nabla I\rangle}.
\end{eqnarray}

We refer to equation~\eqref{eq:crons_ode} as the C-RONS equation. Note that, unlike S-RONS, forming the C-RONS equation does not require computing any functional inner products; it only needs point-wise evaluation at the collocation points $\vc x_i$. Therefore, C-RONS is more computationally efficient. However, this lower computational cost comes at the expense of accuracy since only the residual error at the collocation points is minimized, whereas S-RONS minimizes the error over the entire state space.

We conclude this section by commenting on the special case of S-RONS where the Hilbert space $H$ is the space of square integrable functions $L^2(\mathbb R^d)$ and no regularization is used ($\alpha=0$). Assume that we approximate the inner products~\eqref{eq:srons_M} and~\eqref{eq:srons_f} using Monte Carlo integration instead of symbolic computing, where the Monte Carlo samples $\vc x_i\in\mathbb R^d$ are drawn at random. In this case, the S-RONS equation~\eqref{eq:srons} coincides with the C-RONS equation~\eqref{eq:crons_ode}. We refer to~\cite{anderson2023} for the equivalence proof. Du and Zaki~\cite{du21} used this Monte Carlo approximation with the collocation points distributed according to the uniform Lebesgue measure on $\mathbb R^d$. Bruna et al.~\cite{bruna22} proposed an adaptive sampling method where they draw their collocation points from a distribution which evolves in time. We point out that neither~\cite{du21} nor~\cite{bruna22} use regularization or ensure preservation of conserved quantities.

\section{Choice of the Metric}\label{sec:metric}
Recall that in symbolic RONS, the residual error is minimized with respect to the norm $\|\cdot\|_H$ defined on the Hilbert space $H$. So far, we have stated the results for a general Hilbert space. In this section, we determine the specific choice of the Hilbert space for the Fokker--Planck equation. 

\begin{figure}
	\centering
	\includegraphics[width=\textwidth]{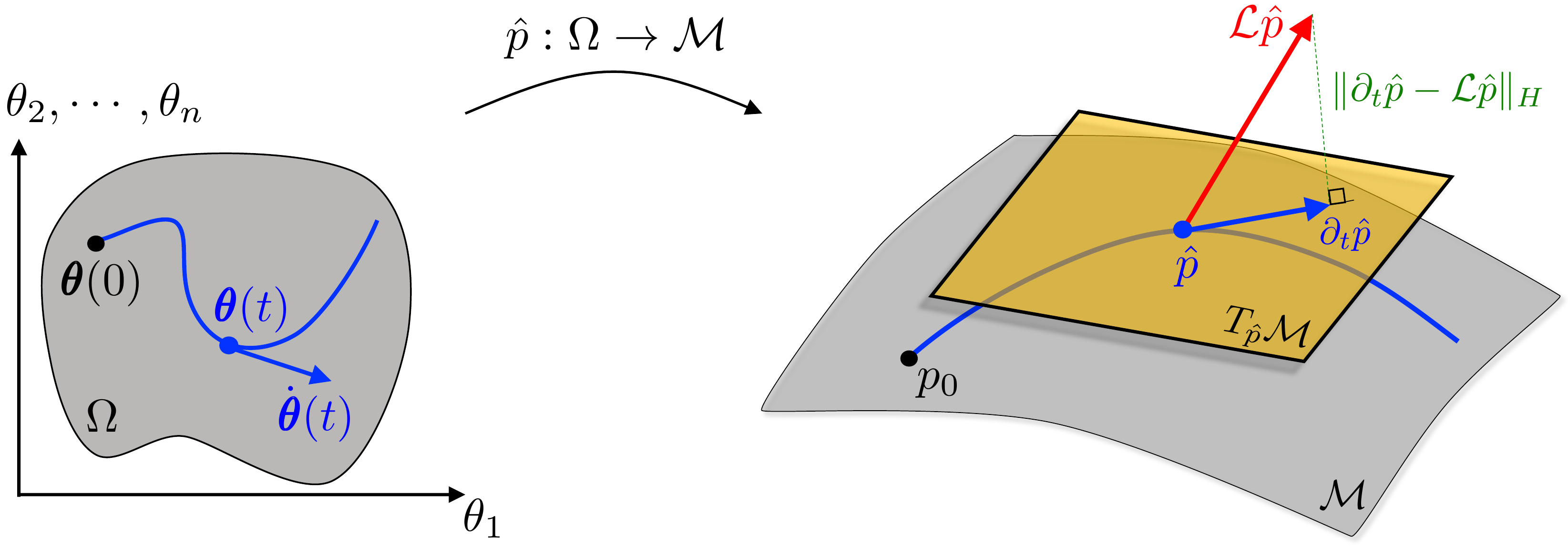}
	\caption{Geometric illustration of the shape-morphing approximate solutions. The image of $\phat$ is a statistical manifold $\mathcal M\subset L^1(\mathbb R^d)$. Any parameter value $\btheta$ defines a point $\phat(\cdot,\btheta)$ on this manifold, with the corresponding tangent space $T_{\phat}\mathcal M$.}
	\label{fig:schem_geom}
\end{figure}

To this end, we first view the approximate solution~\eqref{eq:gen_approx} as a map from the parameters $\btheta\in\mathbb R^n$ to the space of probability densities $L^1(\mathbb R^d)$, 
\begin{align}
\hat p:\ & \Omega \to L^1(\mathbb R^d)\nonumber\\
& \btheta \mapsto \hat p(\cdot,\btheta),
\end{align}
where $\Omega\in\mathbb R^n$ is the set of all admissible parameters $\btheta$. As illustrated in figure~\ref{fig:schem_geom}, the image of the map $\hat p$ forms an $n$-dimensional subset of the infinite-dimensional function space $L^1(\mathbb R^d)$.
In fact, under certain assumptions, the image of $\hat p$ is an immersed manifold $\mathcal M$~\cite{anderson2021}. Note that, although the parameter $\btheta$ is a finite-dimensional vector, the image $\hat p(\cdot,\btheta)$ is a function of $\vc x$, belonging to the infinite-dimensional space $L^1(\mathbb R^d)$.

The set $\mathcal M$ is a statistical manifold in the sense that every point on it is a probability density function~\cite{murray1993}. 
The intrinsic metric associated with a statistical manifold is the so-called Fisher information metric~\cite{Fisher1922}. More specifically, the metric tensor associated with the Fisher information metric is given by 
\begin{equation}
g_{ij}(\btheta) = \int_{\R^d} \frac{\partial \log \phat}{ \partial \theta_i } \frac{\partial \log \phat}{ \partial \theta_j }  \phat\, \id \vc x = \int_{\R^d} \frac{1}{\phat (\vc x,\btheta)} \pard{\phat}{\theta_i}(\vc x,\btheta)\pard{\phat}{\theta_j}(\vc x,\btheta)\id \vc x.
\label{eq:fisher}
\end{equation}
Measuring the distance between two probability distributions $\hat p(\cdot,\btheta_1)$ and $\hat p(\cdot,\btheta_2)$ according to the metric tensor~\eqref{eq:fisher} returns the Fisher--Rao distance between parameterized probability distributions~\cite{Rao1945,Rao1982}.

Therefore, a suitable metric defined on the statistical manifold $\mathcal M$ is the Fisher information metric. Next, we return to the RONS equation~\eqref{eq:srons}. Recall that we stated the results for a
general Hilbert space $H$. Now, let us consider the specific Hilbert space $H=L^2_\mu(\mathbb R^d)$, where $\mu = \hat p^{-1}\id \vc x$ is a weighted Lebesgue measure. In this case, the metric tensor~\eqref{eq:srons_M} can be written explicitly as
\begin{equation}
M_{ij}(\btheta) = \left\langle \frac{\partial \phat}{ \partial \theta_i }, \frac{\partial \phat}{ \partial \theta_j } \right\rangle_{L^2_\mu} = \int_{ \R^d} \frac{1}{\phat} \frac{\partial \phat}{ \partial \theta_i } \frac{\partial \phat}{ \partial \theta_j } \ \id \vc x.
\label{eq:M_wL2}
\end{equation}
We notice that using the Hilbert space $H=L^2_\mu(\mathbb R^d)$, the RONS metric tensor $M_{ij}$ coincides with Fisher information metric~\eqref{eq:fisher}. In other words, taking $H=L^2_\mu(\mathbb R^d)$ induces the Fisher information matrix on the manifold $\mathcal M$ defined by the RONS approximate solution $\hat p$. Therefore, for the Fokker--Planck equation, a suitable choice of the Hilbert space is the weighted Lebesgue space $L^2_\mu(\mathbb R^d)$.

We point out that Bruna et al.~\cite{bruna22} had already proposed this weighted Lebesgue space in an ad hoc manner in their adaptive Monte Carlo estimation of the integrals~\eqref{eq:srons_M}. 
It is interesting that this adaptive sampling can be rigorously justified in the case of Fokker--Planck equation.

In our experience, using symbolic computation to obtain a closed-form expression for~\eqref{eq:M_wL2} is not always feasible. In such cases, we use the Hilbert space $H=L^2(\mathbb R^d)$ to ensure symbolic computing is feasible at the cost of sacrificing the connection between RONS and the Fisher information metric. In Section~\ref{subsec:1D}, we discuss the ramification of this trade-off on a specific example.

\section{Numerical Results} \label{sec:results}
In this section, we assess the accuracy and computational cost of our method on a number of SDEs with progressively higher level of complexity. In all cases, we use Mathematica for symbolic computing and Matlab for numerical time integration of the RONS equations.

\subsection{Benchmark example: Ornstein--Uhlenbeck process}
We consider a one-dimensional (1D) Ornstein--Uhlenbeck (OU) process and show that a Gaussian evolved according to RONS coincides exactly with the true solution of the Fokker--Planck equation corresponding to the OU process.

The Ornstein-Uhlenbeck process $X(t)$ satisfies the SDE,
\begin{equation}
\id X = -\gamma X \id t + \sigma \id W, \quad X(0) = 0,
\label{eq:OU}
\end{equation}
where $\gamma>0$ is the drift coefficient and $\sigma >0$ is noise intensity.
The Fokker--Planck equation associated with the OU process~\eqref{eq:OU} is
\begin{equation}
\frac{ \partial p }{\partial t} = \gamma \frac{ \partial }{ \partial x } \big( x p \big) + \frac{ \sigma^2 }{ 2 } \frac{ \partial^2 p }{ \partial x^2 },\quad p\left(x,0\right)=\delta\left(x\right),
\label{eq:OU_FPE}
\end{equation}
where $\delta(x)$ is the Dirac delta function centered at the origin.
The Fokker--Planck equation~\eqref{eq:OU_FPE} admits the exact solution~\cite{hanggi2007},
\begin{equation}
p(x,t) = \sqrt{ \frac{ \gamma }{ \pi \sigma^2 ( 1 - \exp[ -2 \gamma t ] ) } } \exp \bigg[  - \frac{ \gamma x^2 }{ \sigma^2 ( 1 - \exp[ -2 \gamma t ] ) } \bigg],
\label{eq:OU_FPE_soln}
\end{equation}
which is a Gaussian whose amplitude decays over time while its variance grows. Also note that $p(x,t)$ tends to the Dirac delta function $\delta (x)$ as time $t$ tends to zero.

To apply RONS, we consider the Gaussian solution,
\begin{equation}
\hat{p} (x, \pmb \theta(t)) = A(t) \exp \bigg[ -\frac{ (x - c(t))^2 }{ L^2(t) } \bigg],
\label{eq:gauss_ansatz_1D}
\end{equation}
with the time-dependent parameters $\pmb \theta (t) = (A(t), L(t), c(t))^\top$. 
Note that this is the Gaussian approximate solution~\eqref{eq:gen_approx} with only one mode $(r = 1)$. Here, we do not square the amplitude to simplify the following analysis. 

As mentioned in Section \ref{sec:prelim}, the approximate solution $\hat  p$ is a probability density function (PDF), and so we must ensure that total probability of our approximate solution is always equal to one. 
This is a conserved quantity of the Fokker--Planck equation which we enforce when applying RONS by ensuring
\begin{equation}
I(\pmb \theta(t)) = \int_{ \R } \hat p(x, \pmb \theta(t)) \ \id x = \sqrt{\pi} A(t) L(t) = 1, \quad \forall t\geq 0.
\label{eq:OU_const_dens}
\end{equation}

Applying RONS to the Fokker--Planck equation associated with the OU process, using the Gaussian approximate solution~\eqref{eq:gauss_ansatz_1D} and the Hilbert space  $H = L^2(\R)$,
the corresponding S-RONS equation~\eqref{eq:srons} reads
\begin{equation}
\dot{  A }= A \bigg( \gamma - \frac{ \sigma^2}{ L^2 }\bigg) , \quad \dot{  L }=  \frac{ \sigma^2}{ L } - \gamma L, \quad \dot{  c } = -\gamma c.
\label{eq:OU_ODEs}
\end{equation}

To solve these ODEs, we need to specify the appropriate initial conditions. For the Gaussian~\eqref{eq:gauss_ansatz_1D} to coincide with the initial condition of the Fokker--Planck equation~\eqref{eq:OU_FPE} at time $t=0$, we choose the initial parameter values,
\begin{equation}
A(t_0) = \sqrt{ \frac{ \gamma }{ \pi \sigma^2 ( 1 - \exp[ -2 \gamma t_0 ] ) } }, \quad L(t_0) = \frac{1}{\sqrt{\pi}A(t_0)}, \quad c(t_0) = 0.
\label{eq:paramvals}
\end{equation}
Note that as $t_0\to 0$, the initial condition $\phat (\vc x,\btheta(t_0))$ approaches the Dirac delta function $\delta (x)$ as required. 
The exact solution to the S-RONS equation~\eqref{eq:OU_ODEs}, with the initial condition~\eqref{eq:paramvals}, is given by
\begin{equation}
A(t) = \sqrt{ \frac{ \gamma }{ \pi \sigma^2 ( 1 - \exp[ -2 \gamma t ] ) } }, \quad L (t)= \bigg( \frac{ \gamma }{  \sigma^2 ( 1 - \exp[ -2 \gamma t ] ) } \bigg)^{-1/2} , \quad c(t) = 0.
\label{eq:OU_params}
\end{equation}

Substituting this solution into $\phat (\vc x,\btheta(t))$, we recover the exact solution~\eqref{eq:OU_FPE_soln} to the Fokker--Planck equation.
This benchmark example shows that the RONS solution coincides with the exact solution of the Fokker--Planck equation for the OU process.

\subsection{One-dimensional bistable potential}
\label{subsec:1D}

In this section, we consider an SDE where the dynamics are driven by a potential function. Our main focus here is to highlight the impact of the choice of the Hilbert space $H$ and the regularization parameter $\alpha$. We consider the SDE
\begin{equation}
\id  X(t)  = - V^\prime (x) \ \id t + \sigma \id W,
\label{eq:bistable_sde}
\end{equation}
where $V: \R \to \R$ is the potential.
The Fokker--Planck equation corresponding to~\eqref{eq:bistable_sde} reads
\begin{equation}
\frac{ \partial p }{ \partial t } = \frac{ \partial  }{ \partial x } \big( V^\prime(x) p \big) + \nu \frac{ \partial^2 p }{ \partial x^2 }.
\label{eq:bistable_FPE}
\end{equation}
In general, the analytical solution for this Fokker-Plank equation for all times is not known.
However, the asymptotically stable steady state solution is given by
\begin{equation}
p_{\text{eq}}(x) = C \exp \bigg[ - \frac{ V(x) }{ \nu } \bigg] ,
\label{eq:bistable_equil}
\end{equation}
where $C$ is a normalizing constant.

Here, we consider the potential  
\begin{equation}
V(x) =  \frac{x^4}{4}  - \frac{x^2}{2}.
\label{eq:bistable_V}
\end{equation} 
This potential is symmetric with two minima at $x = \pm 1$, and so the equilibrium solution~\eqref{eq:bistable_equil} is bimodal with peaks at $x = \pm 1$. Similar bistable potentials have been studied by several other authors~\cite{bernstein1984, farazmand2020,Mamis2021a,wei1999}.

To apply RONS we once again consider the Gaussian ansatz~\eqref{eq:gen_approx}, where we will now use sums of Gaussians rather than a single Gaussian in our approximation. 
As before, we also enforce that total probability of the approximate solution is always equal to 1 to ensure that $\hat{p}$ is in fact a probability density function.

In addition to exploring the effects from changing the number of modes used in our approximate solution, we also study the the choice of Hilbert space $H$ for our inner products.
In particular, we compare the results for the Hilbert space  $H = L^2(\R)$ and the weighted Hilbert space  $H = L^2_{\mu}(\R)$. As discussed in Section~\ref{sec:metric},
when using the weighted inner product, the metric tensor $M$ coincides with the Fisher information matrix.
\begin{figure}
	\includegraphics[width=\textwidth]{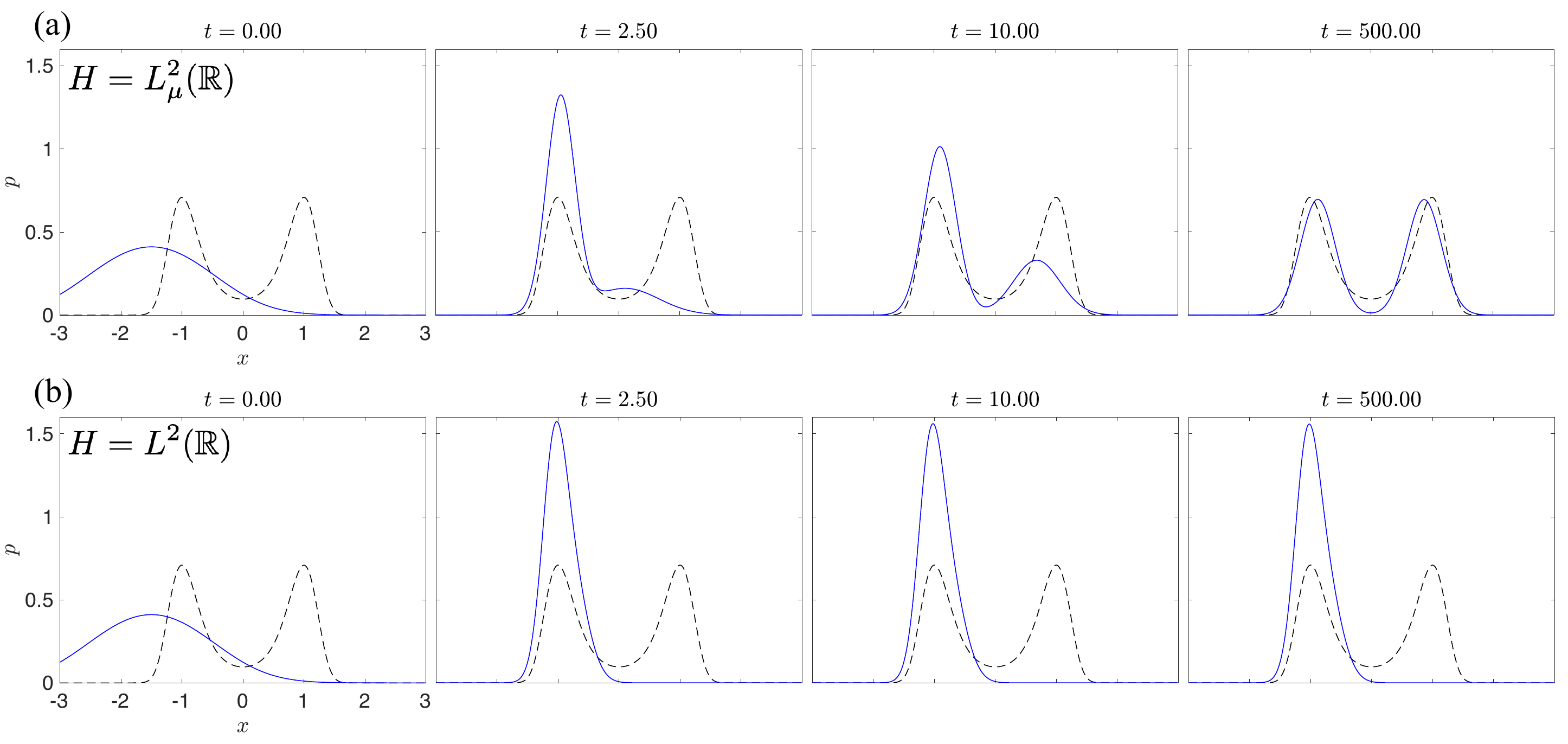}
	\caption{Evolution of $\hat{p}(x, \pmb \theta (t))$ from applying RONS using 2 Gaussians in the approximate solution (blue curves). The true equilibrium density $p_{\text{eq}}(x)$ is marked by the dashed black curve. The initial condition is $A_i(0) = 1/2$, $L_i(0) = 2/ \sqrt\pi$, $c_1(0) = -1$, $c_2(0) =-2$.  The Hilbert space is (a) $H=L^2_{\mu}(\R)$  and (b) $H=L^2(\R)$.}
	\label{fig:1D_2gauss}
\end{figure}

Using the Hilbert space $H=L^2(\R)$, we are able to symbolically calculate the inner products of the RONS equations and apply S-RONS.
This approach is scalable and allows for rapid time integration of the ODEs.
In contrast, when using the Hilbert space $H=L^2_{\mu}(\R)$, obtaining closed-form symbolic expressions for the inner products was not possible. Consequently, we resort to using C-RONS for the weighted $L^2$ inner product space. Note that C-RONS requires sampling which can be expensive in higher dimensions, but it is not an obstacle in this 1D problem.

For this section, we integrate the RONS equation using Matlab's \texttt{ode15s} solver \cite{shampine1997}.
In our numerical experiments, \texttt{ode15s} takes large time steps once the approximate solution $\hat p$ is near its equilibrium solution $\hat p_{\text{eq}}$.
This allows for rapid simulations over long time scales, which helps us obtain the equilibrium solution predicted by RONS.

We first apply RONS to the Fokker--Planck equation~\eqref{eq:bistable_FPE} with the bistable potential, using $r=2$ modes for the Gaussian approximate solution~\eqref{eq:gen_approx}.
We use a noise intensity of $\sigma = 0.5$ for all simulations in this section.
Figure~\ref{fig:1D_2gauss} compares the evolution of the approximate solution $\hat p$ predicted by S-RONS using the $L^2$ inner product and C-RONS using the weighted $L^2_{\mu}$ inner product.
When applying C-RONS we take 100 equidistant collocation points on the interval $x = [-4,4]$.

\begin{figure}
	\includegraphics[width=\textwidth]{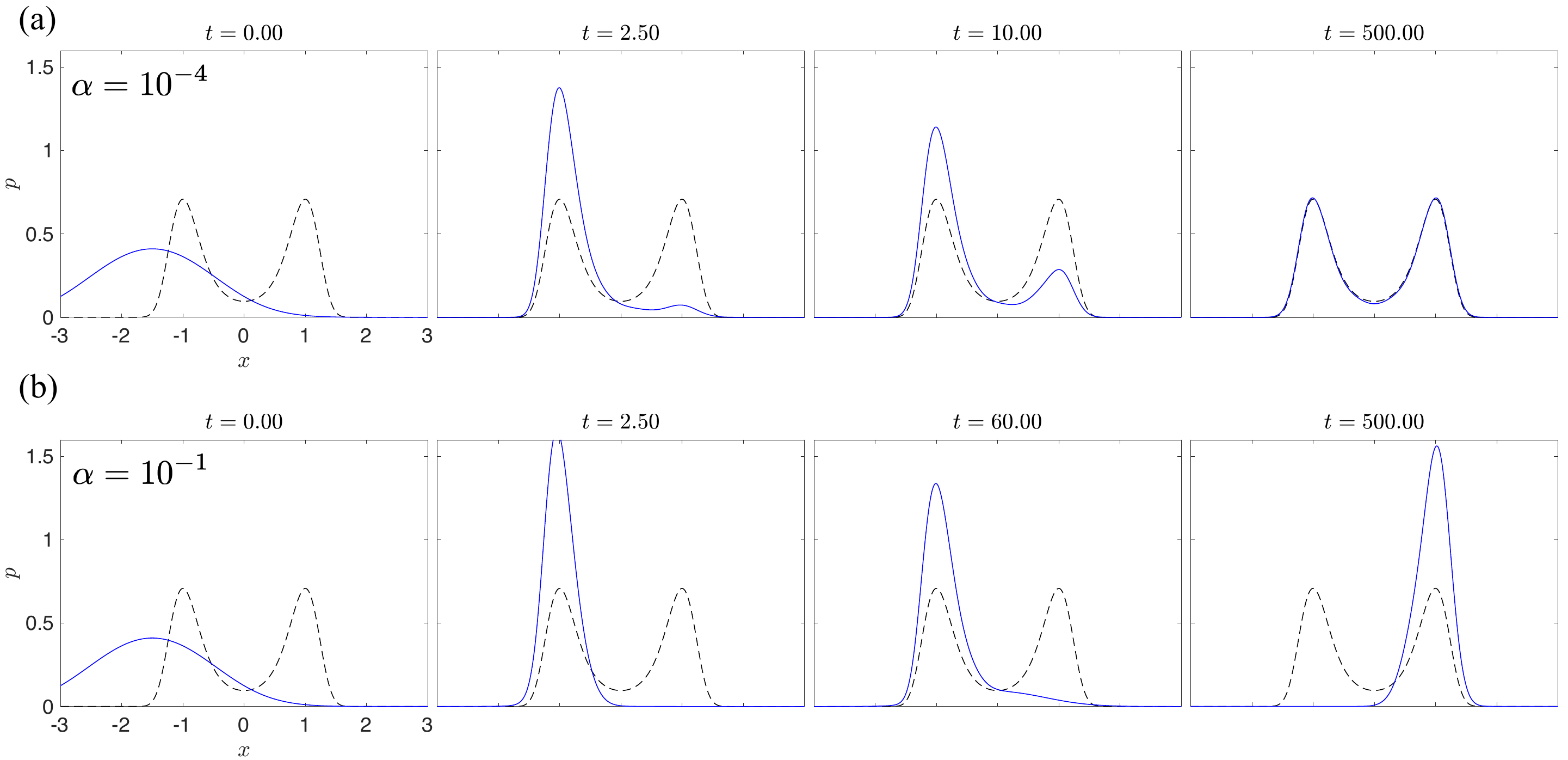}
	\caption{Evolution of $\hat{p}(x, \pmb \theta (t))$ from applying RONS with the Hilbert space $H = L^2(\R)$ and using 10 Gaussians in the approximate solution (blue curves). We also mark the true equilibrium density $p_{\text{eq}}(x)$ with a dashed black curve. The initial condition is $A_i(0) = 1/\sqrt{20}$, $L_i(0) = 2/\sqrt \pi$, where the amplitudes are chosen so that total probability is 1. Half of the Gaussians are initially placed at $x = -1$, with the other half placed at $x =-2$.  Regularization parameter is (a) $\alpha = 10^{-4}$ and (b) $\alpha = 10^{-1}$.}
	\label{fig:1D_10gauss}
\end{figure}

First we observe that using the Hilbert space $H = L^2_{\mu}(\R)$ produces an approximate solution which is a reasonable approximation of the equilibrium density by two Gaussians. 
However, taking $H = L^2(\R)$ fails in this case as both Gaussians converge to the same peak at $x=-1$.
We have observed the same behavior for a wide range of initial conditions for the two Gaussians.
When using $H = L^2_{\mu}(\R)$, corresponding to the Fisher information metric, the approximate solution always converges to the true equilibrium density, 
whereas using $H = L^2(\R)$ with two Gaussians leads to the incorrect equilibrium.
The only exception to this is when we start from initial conditions which are symmetric, with the Gaussians initially placed on the opposite sides of the origin.
In this particular case, the approximate solution converges to a reasonable approximation of the true equilibrium even when we use $H = L^2(\R)$ (not shown here).

The fact that using the weighted inner product space $H = L^2_{\mu}(\R)$ leads to better results is not surprising. As discussed in Section~\ref{sec:metric}, this is equivalent to using the Fisher information metric on the manifold of the approximate solution which is the natural metric for a statistical manifold.
Nonetheless, using the unweighted Hilbert space $L^2(\R)$ is still desirable since it allows us to use S-RONS which requires no sampling and incurs no numerical error in approximating the inner products.
So, can we somehow fix the issue of converging to the wrong equilibrium solution using $H=L^2(\R)$? The answer is yes as long as the number of terms $r$ in the approximate solution~\eqref{eq:gen_approx} is large enough.

For instance, let us consider the same initial condition used to produce figure~\ref{fig:1D_2gauss}, but we now use $r=10$ Gaussians in our approximate solution.
As discussed in Section~\ref{sec:evolution}, when using a large number of parameters in the approximate solution, we must apply Tikhonov regularization to alleviate the stiffness of the RONS ODEs.
Figure~\ref{fig:1D_10gauss}(a) shows the evolution of the approximate solution $\hat{p}$ with $H = L^2(\R)$ and the regularization parameter $\alpha = 10^{-4}$. Unlike the previous case where only two Gaussians were used, the approximate solution converges to the correct equilibrium density. Of course, one should be cautious not to over-regularize the problem. As shown in figure~\ref{fig:1D_10gauss}(b), choosing $\alpha = 10^{-1}$ causes the solution to approach the incorrect equilibrium again. Although here we chose the Tikhonov regularization parameter $\alpha=10^{-4}$ by trial and error, there exist more rigorous methods for choosing this parameter a priori~\cite{engl1996,Ito2011}.
\begin{figure}
	\centering
	\includegraphics[width=0.4\textwidth]{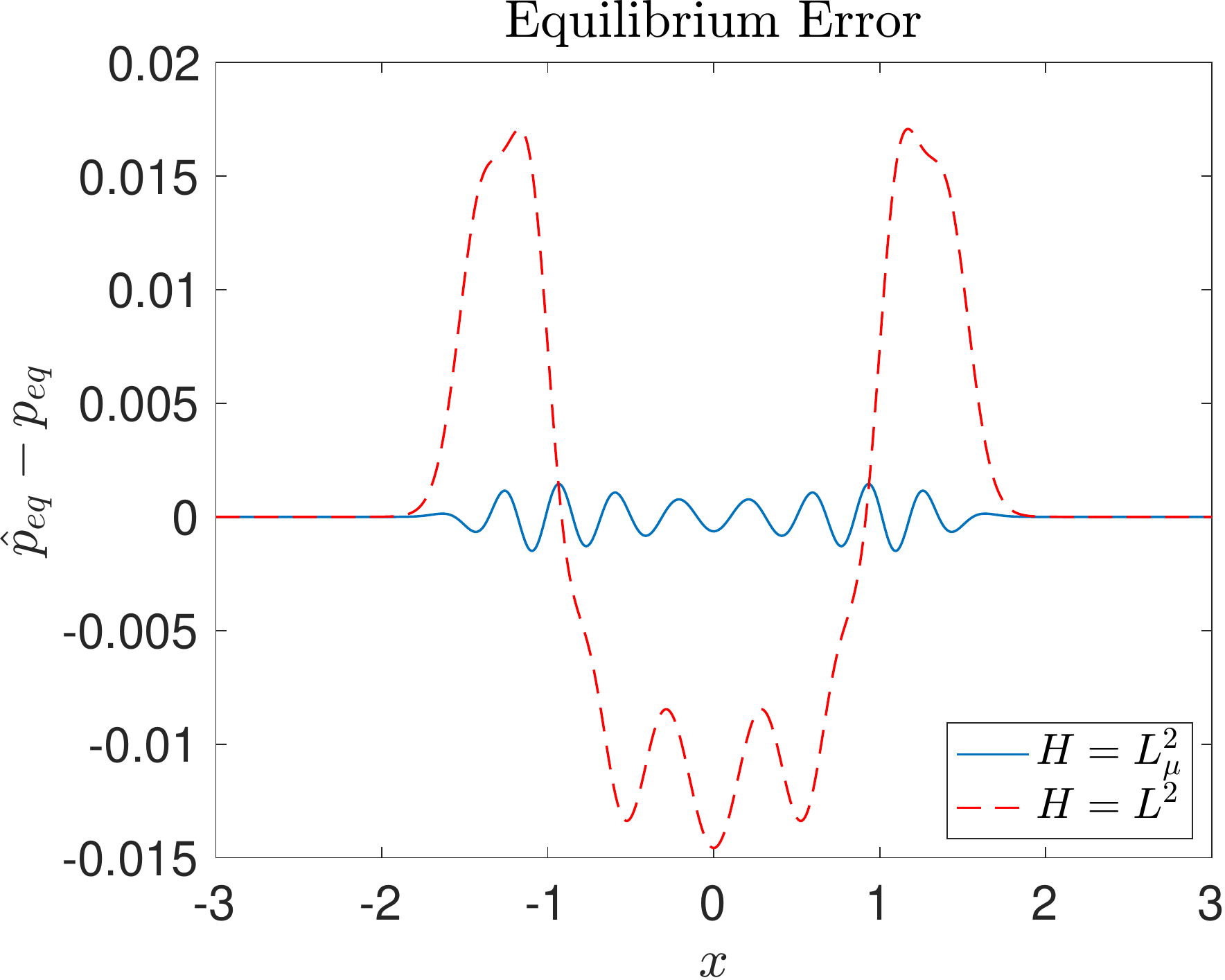}
	\caption{Error in approximating the true equilibrium density $p_{\text{eq}}$ using 10 Gaussians in the approximate solution, regularization parameter $\alpha = 10^{-4}$, and two different choices of Hilbert space. The initial conditions are  $A_i(0) = 1/\sqrt{20}$, $L_i(0) = 2/\sqrt \pi$.}
	\label{fig:1D_10gauss_error}
\end{figure}

Figure~\ref{fig:1D_10gauss_error} shows the equilibrium error when using $r=10$ Gaussians with both the weighted and unweighted inner products.
A comparison of computation time and errors is also provided in Table \ref{table:1d_bistable_comptime}.
Although in both cases the errors are small, using the Hilbert space $L^2_\mu(\R)$ provides a more accurate solution than the Hilbert space $L^2(\R)$.
However, as we noted earlier, using the unweighted inner product space $L^2(\R)$ is scalable to higher dimensions since it allows the use of symbolic computing instead of collocation points (cf. Section~\ref{sec:harmonic} below).
\begin{table}
	\centering
	\begin{tabular}{|l||ccc|}
		\hline
		\textbf{Bistable Potential $(r = 10)$} &
		\begin{tabular}[c]{@{}c@{}}Symbolic \\ computation\end{tabular} &
		\begin{tabular}[c]{@{}c@{}}Time \\ integration\end{tabular}  &
		\begin{tabular}[c]{@{}c@{}}Relative error \\ of equilibrium\end{tabular}\\ \hline
		\multicolumn{1}{|l||}{C-RONS ($H = L^2_\mu (\R)$)} &
		none &
		2.11 seconds &
		0.002 \\ \hline
		\multicolumn{1}{|l||}{S-RONS ($H = L^2 (\R)$)} &
		2.78 seconds &
		1.50 seconds &
		0.03 \\
		\hline
	\end{tabular}
	\caption{Comparison of computational time and errors for the 1D bistable potential using symbolic RONS (S-RONS) and collocation RONS (C-RONS).}
	\label{table:1d_bistable_comptime}
\end{table}

\subsection{Stochastic Duffing oscillator}

In this section, we consider the stochastic Duffing oscillator~\cite{kumar2006,pradlwarter2001,spencer1993} excited by white noise,
\begin{equation}
\id \vc X = \begin{pmatrix}
y\\
a_1 x + a_2 y + a_3 x^3
\end{pmatrix} \id t + \sigma \begin{pmatrix}
0 & 0 \\
0 & 1
\end{pmatrix} \id \vc W.
\label{eq:Duffing_SDE}
\end{equation}
Here $\vc X = (x, y )^\top$ where $x(t)$ is the displacement and $y(t)=\dot x(t)$ is the velocity. The vector $\vc  W(t) = (W_1(t), W_2(t))^\top$ represents the standard Wiener process in two dimensions. 
The coefficients $a_i$ are constants where $a_1$ controls the stiffness of the oscillator, $a_2$ controls the damping, and $a_3$ controls the strength of the nonlinearity in the restoring force of the oscillator. The Fokker--Planck equation for the stochastic Duffing oscillator is given by
\begin{align}
\pard{p}{t} &= - \bigg[  y \pard{p}{x} + a_2 p + ( a_1 x + a_2 y + a_3 x^3 ) \pard{p}{y} \bigg] + \frac{ \sigma^2 }{ 2 } \pard{^2p}{y^2}.
\label{eq:Duffing_FP}
\end{align}

An analytic solution to the Fokker--Planck equation \eqref{eq:Duffing_FP} is not known for all times.
However, the asymptotically stable equilibrium solution is given by~\cite{pradlwarter2001}
\begin{equation}
p_{ \text{eq} }(\vc x) = C \exp \bigg[ \frac{ -a_1 a_2 x^2 - \frac{ 1 }{ 2 } a_2 a_3 x^4 + a_2 y^2 }{ \sigma^2 } \bigg],
\label{eq:Duffing_equil}
\end{equation}
where $C$ is a normalizing constant.
Following \cite{pradlwarter2001}, we use the parameter values $( a_1, a_2, a_3 ) = ( 1, -0.2, -1 )$ and a noise intensity of $\sigma = 1/\sqrt{20}$.
This leads to a bimodal equilibrium distribution, with peaks at $(x,y) = (\pm 1, 0)$.

To approximate solutions of the Fokker--Planck equation~\eqref{eq:Duffing_FP}, we use the Hilbert space $H = L^2(\R^2)$ and the Gaussian approximate solution~\eqref{eq:gen_approx}.
We evolve the parameters using the S-RONS equation~\eqref{eq:srons} with the regularization parameter $\alpha = 10^{-3}$. These ODEs are integrated numerically using Matlab's \texttt{ode45} \cite{dormand1980,shampine1997}.

\begin{figure}
	\centering
	\includegraphics[width=\textwidth]{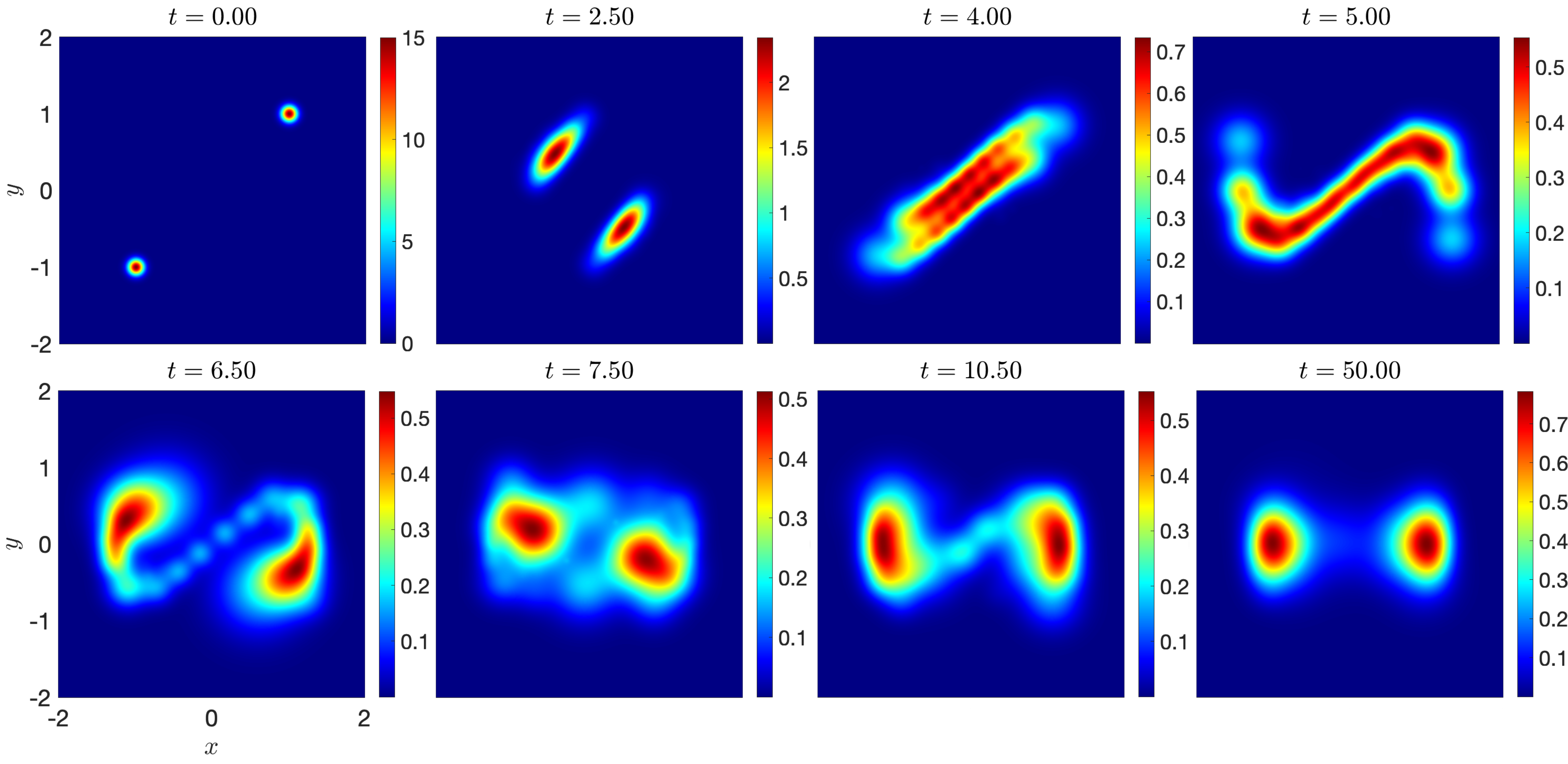}	
	\caption{PDF predicted by applying RONS to the Fokker--Planck equation~\eqref{eq:Duffing_FP} for the Duffing oscillator excited by white noise. We use 30 Gaussians in the approximate solution. Initial conditions are $A_i(0) = 1$, $L_i(0) = ( 30\pi )^{-1/2}$, and half the Gaussians are placed at $(x,y) = (-1,-1)$ with the other half at $(x,y) = (1,1) $. }
	\label{fig:Duffing_snapshots_RONS}
\end{figure}

\begin{figure}
	\centering
	\includegraphics[width=\textwidth]{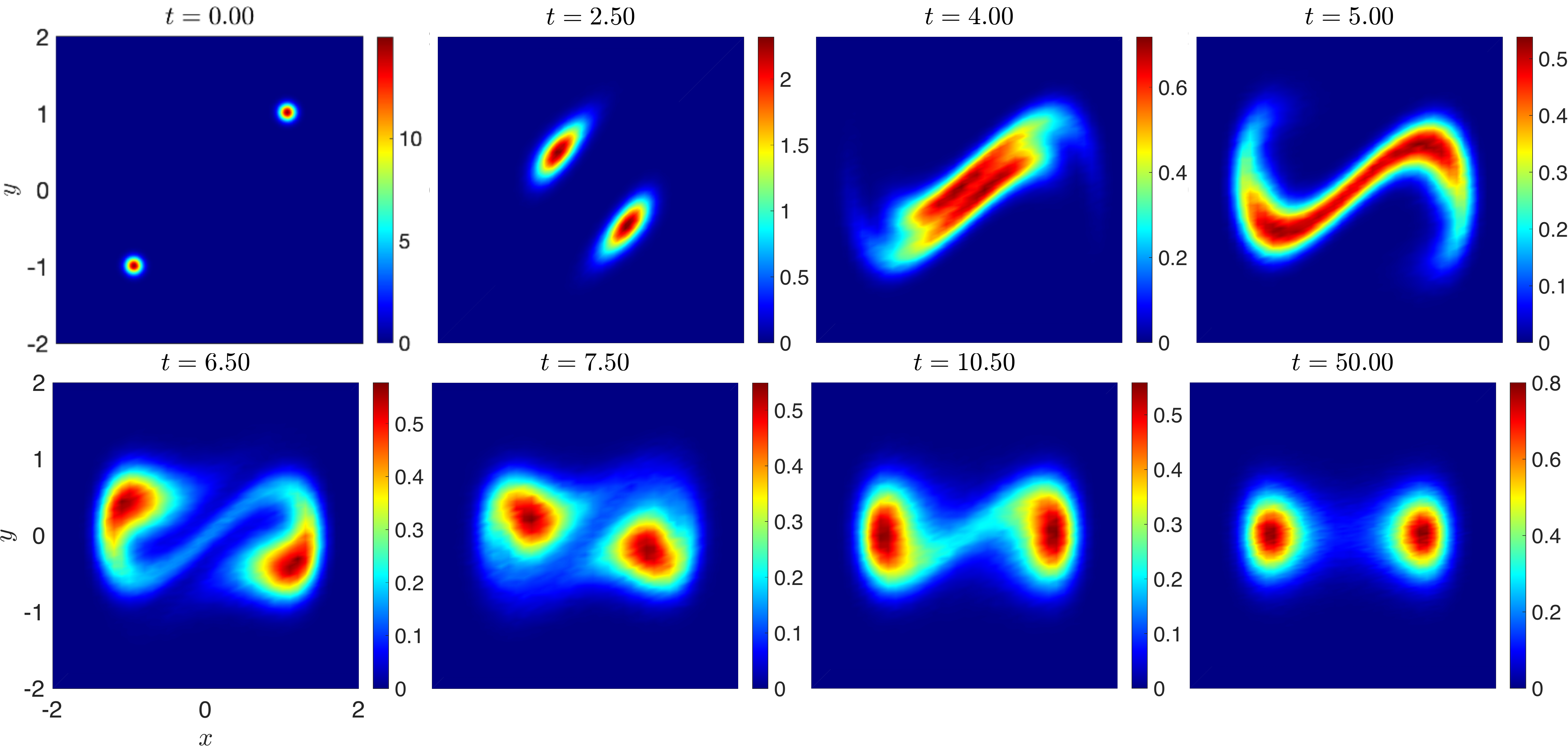}
	\caption{PDF predicted by direct Monte Carlo simulations of the Duffing oscillator excited by white noise~\eqref{eq:Duffing_SDE}. The initial distribution is the sum of 30 Gaussians with parameter values $A_i(0) = 1$, $L_i(0) = ( 30\pi )^{-1/2}$, and half the Gaussians are placed at $(x,y) = (-1,-1)$ with the other half at $(x,y) = (1,1) $. }
	\label{fig:Duffing_snapshots_MC}
\end{figure}

Figure \ref{fig:Duffing_snapshots_RONS} shows the evolution of the approximate solution $\hat p(\vc x, \pmb \theta (t))$ with $r=30$ modes. 
The initial condition is set up so that $15$ Gaussians are centered at $(x,y)=(-1,-1)$ and the other $15$ are centered at $(x,y)=(+1,+1)$, leading to a bimodal initial condition.
The approximate solution is evolved until it converges to the equilibrium density and virtually no further change is detected.

We compare the S-RONS solution against large-scale Monte Carlo simulations of the Duffing SDE~\eqref{eq:Duffing_SDE}. 
To this end, we evolve $10^6$ particles using the predictor-corrector scheme of Ref.~\cite{cao2015}. The particles are drawn at random such that their distribution matches that of the S-RONS simulations at the initial time $t=0$. The evolution of the resulting Monte Carlo PDF is shown in figure~\ref{fig:Duffing_snapshots_MC}.

Comparing figures~\ref{fig:Duffing_snapshots_RONS} and~\ref{fig:Duffing_snapshots_MC}, we observe that S-RONS, not only returns the correct equilibrium density, but also reproduces the transient dynamics very well. RONS does not perfectly capture some of the finer features seen in the Monte Carlo approach, such as the tails of the solution at $t = 5$. This is expected due to the fact that we are evolving only 30 Gaussians. In fact, increasing the number of terms to $r=100$ allows us to capture these fine features as well (not shown here for brevity).

In terms of computational cost, RONS is over 600 times faster than direct Monte Carlo simulations. As reported in Table \ref{table:Duffing_comptime}, the Monte Carlo simulations take over 2 hours to complete, whereas RONS takes 3.83 seconds for symbolic computing and approximately 10 seconds for time integration. We emphasize that the symbolic computation for S-RONS only needs to be carried out once; changing the initial condition or increasing the number of modes $r$ does not require additional symbolic computing.

\begin{figure}
	\centering
	\includegraphics[width=\textwidth]{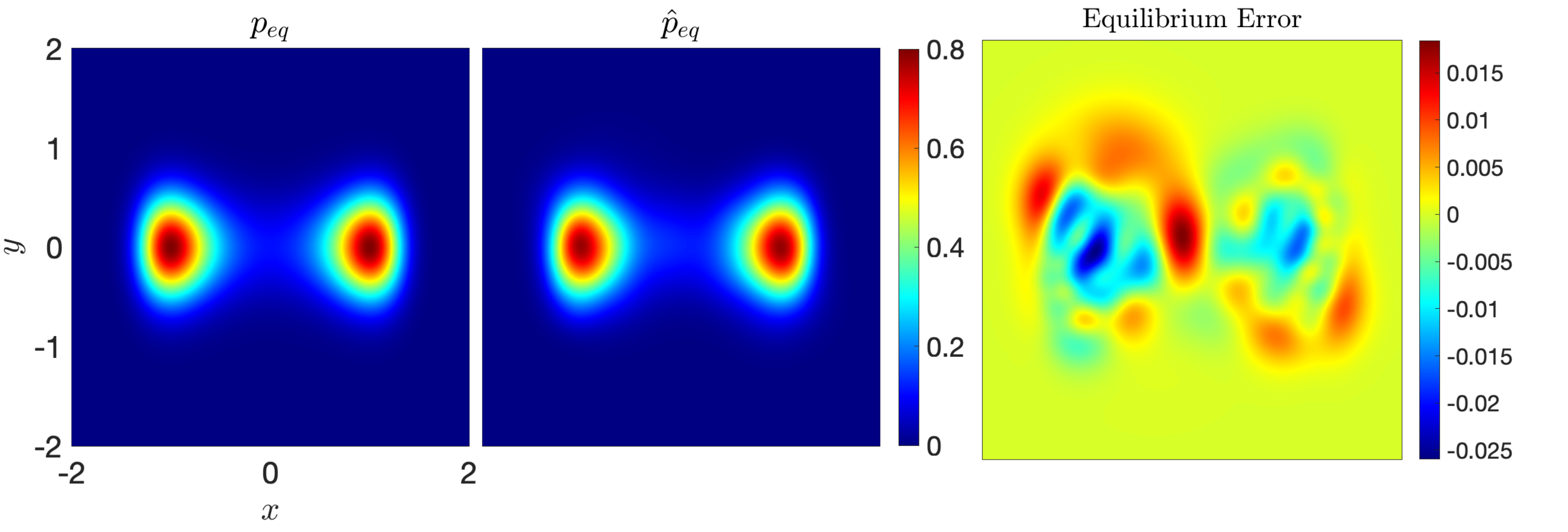}
	\caption{Comparing the true equilibrium density $p_{\text{eq}}$ to the approximate density $\phat_{\text{eq}}$ obtained by RONS using 30 Gaussians. Initial conditions are $A_i(0) = 1$, $L_i(0) = ( 30\pi)^{-1/2}$, with half the Gaussians placed at $(x,y) = (-1,-1)$ and the other half at  $(x,y) = (1,1) $.}
	\label{fig:Duffing_error}
\end{figure}

Finally, figure \ref{fig:Duffing_error} compares the RONS solution at time $t = 50$ with the analytical equilibrium density~\eqref{eq:Duffing_equil}.
We can see that RONS provides an excellent approximation of the analytical equilibrium position, with an $L^2$ relative error of $2.2\%$.

\begin{table}
	\centering
	\begin{tabular}{|l||cc|}
		\hline
		\textbf{Stochastic Duffing Oscillator $(r = 30)$} &
		\begin{tabular}[c]{@{}c@{}}Symbolic \\ computation\end{tabular} &
		\begin{tabular}[c]{@{}c@{}}Time \\ integration\end{tabular}  \\ \hline
		\multicolumn{1}{|l||}{Monte Carlo Simulations} &
		none &
		$ 136.44$ minutes \\ \hline
		\multicolumn{1}{|l||}{Symbolic RONS} &
		3.83 seconds &
		\ \ \ 9.56 seconds \\
		\hline
	\end{tabular}
	\caption{Comparison of computation time for stochastic Duffing oscillator example using symbolic RONS and Monte Carlo simulations of the SDE~\eqref{eq:Duffing_SDE}.}
	\label{table:Duffing_comptime}
\end{table}

\subsection{Harmonic trap}\label{sec:harmonic}
In this section, we study an SDE in eight dimensions driven by a harmonic trap, which was also investigated in \cite{anderson2023,bruna22}.
More specifically, we consider a system of $d$ interacting particles whose motion is governed by the SDE,
\begin{equation}
\id X_i = g(t, X_i) \id t +  \sum_{ j = 1 }^d K(X_i, X_j) \id t +  \sigma \id W_i, \quad i = 1,2,...,d,
\label{eq:harmonic_SDE}
\end{equation}
where $X_i(t)$ denotes the position of the $i$-th particle.
The function $g:[0, \infty]\times \R \to \R$ is a forcing term and $K:\R\times \R \to \R$ describes the interaction between particles.
The corresponding Fokker--Planck equation is given by
\begin{equation}
\frac{\partial p}{ \partial t } = \sum_{ i = 1 }^{d}  - \frac{\partial  }{\partial x_i} \bigg[ \bigg( g(t,x_i)+ \sum_{ j = 1 }^{d} K(x_i,x_j) \bigg) p \bigg]  + \nu \frac{\partial^2 p }{\partial x_i^2},
\label{eq:harmonic_FP}
\end{equation}
where $\nu=\sigma^2/2$.

As in \cite{bruna22}, we set
\begin{equation}
g(t, x_i) = a(t) - x_i, \quad K(x_i,x_j) = \frac{\gamma}{d} (x_j - x_i).
\label{eq:FP_gandk}
\end{equation}
The choice of $g$ corresponds to particles in a harmonic trap centered around $a(t)$ while the particles also attract each other due to interaction term $K$.
A significant advantage of these choices for $g$ and $K$ is that we can directly compute the mean and covariance of the particles to serve as a benchmark for our RONS results.
By taking the expected value of the SDE~\eqref{eq:harmonic_SDE}, we obtain an expression for the mean of each particle
\begin{equation}
\dot{\bar{X}}_i = a(t) - \bar{X}_i + \frac{\gamma}{d}  \sum_{ j = 1 }^d (\bar{X}_j - \bar{X}_i),\quad i = 1,2,...,d,
\label{eq:trap_mean}
\end{equation}
where $\bar{X}_i = \mathbb{E}[X_i]$.
We can similarly derive an expression for the evolution of the correlation matrix $\Sigma_{ij} = \mathbb{E}[X_i X_j]$,
\begin{align}
\dot{\Sigma}_{ij} &= a(t) ( \bar{X}_j + \bar{X}_i ) - 2(1+\gamma)\Sigma_{ij} + \frac{\gamma}{d} \sum_{ l = 1 }^d  ( \Sigma_{lj} + \Sigma_{li} ) + 2  \nu \delta_{ij},& \quad &i,j \in \{1,2,...,d\},
\label{eq:trap_sigma}
\end{align}
where $\delta_{ij}$ denotes the Kronecker delta. 
The covariance matrix, whose entries are given by $\Sigma_{ij} - \bar X_i \bar X_j$, is then calculated using the solutions of~\eqref{eq:trap_mean} and~\eqref{eq:trap_sigma}.

For the Fokker--Planck equation, we choose the initial condition,
\begin{equation}
p(\vc x, 0) = \left(2 \pi \sigma_0^2\right)^{-d/2} \exp \bigg[ -\frac{|\vc x - \pmb \mu|^2}{2\sigma_0^2} \bigg].
\label{eq:harmonic_IC}
\end{equation}
where $\sigma_0^2=0.1$ and $\pmb \mu \in \R^d $ is given by $\mu_i = i - 1$ for $i = 1,...,d,$. 
The remaining parameters are given by $a(t) = 1.25(\sin(\pi t) + 1.5) $, $\gamma = 0.25$, $d = 8$, and $\nu = 0.01$. 

We use S-RONS with the Hilbert space $H=L^2(\R^8)$ to evolve the approximate solution~\eqref{eq:gen_approx}. For this approximate solution to coincide with the initial condition~\eqref{eq:harmonic_IC} at time $t=0$, we choose the parameter values $A_i(0) = (2 \pi \times 0.1)^{-4}r^{-1}$, $L_i^2(0) = 0.2$, and $c_i(0) = \pmb \mu$.
We again enforce the total probability of the approximate solution $\hat{p}(\vc x, \pmb \theta)$ to be always equal to one.
\begin{figure}
	\centering
	\includegraphics[width=0.49\textwidth]{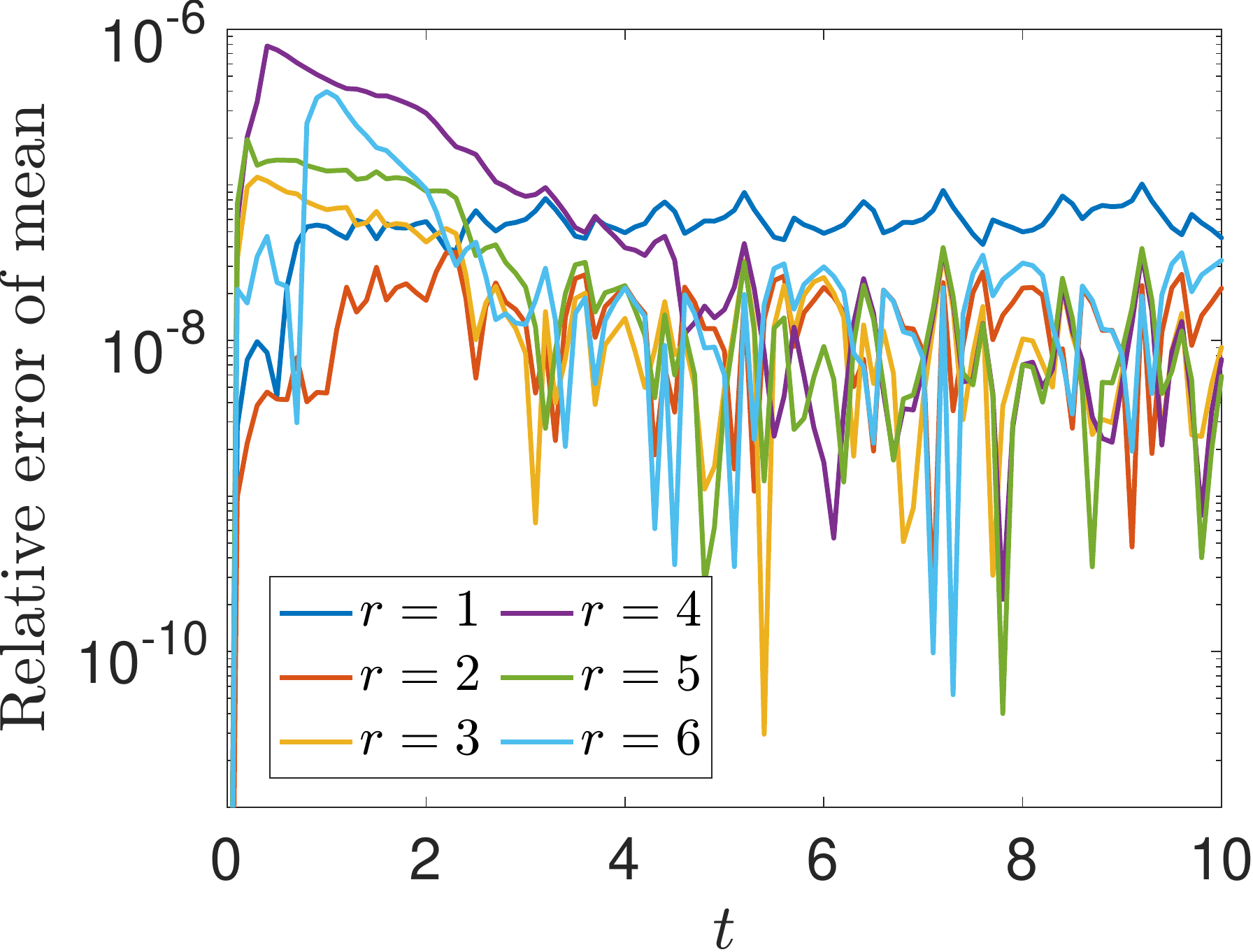}
	\includegraphics[width=0.49\textwidth]{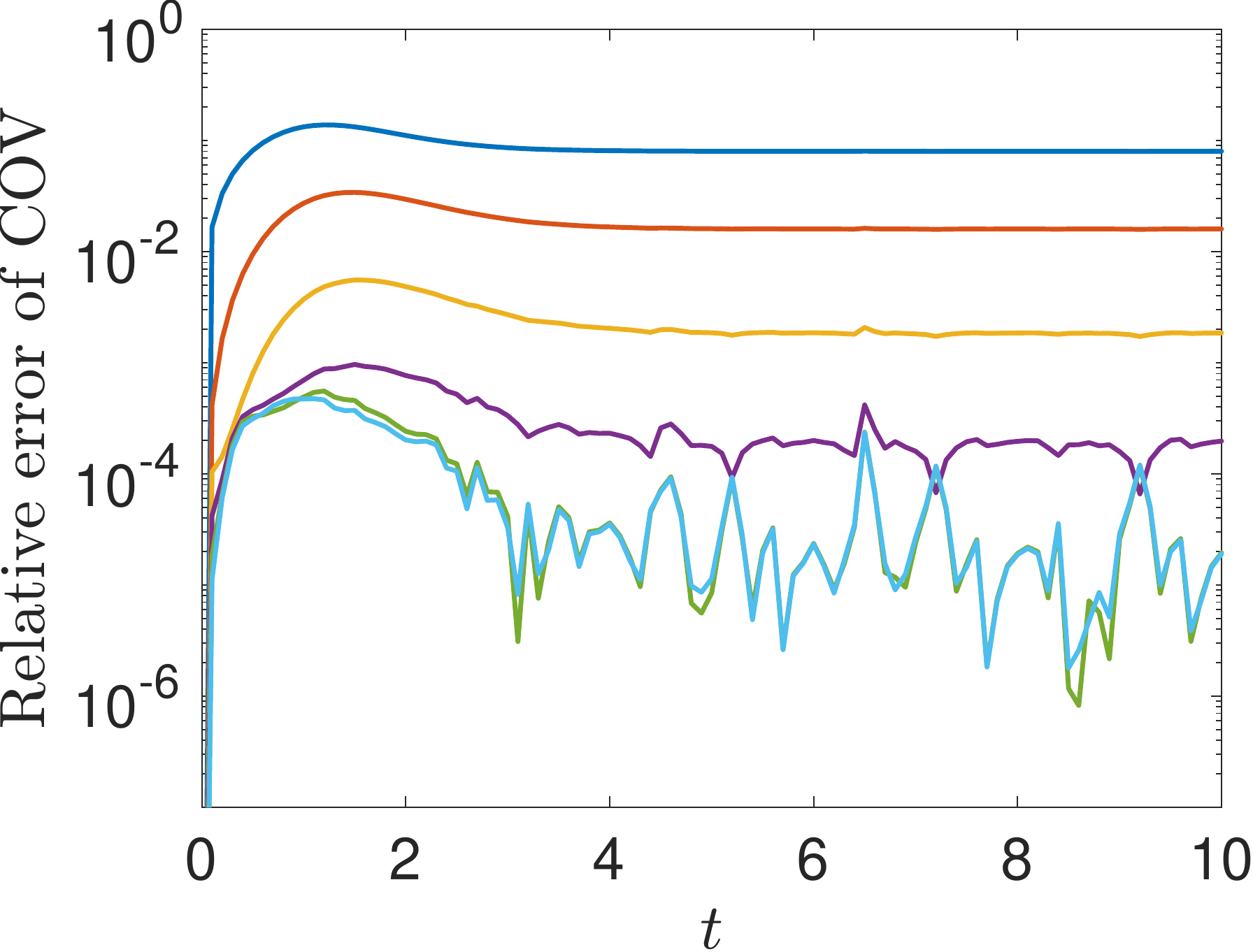}
	\caption{RONS simulation of the harmonic trap using various numbers of modes with regularization parameter $\alpha = 10^{-8}$. }
	\label{fig:harmonic_error_chap3}
\end{figure}

We numerically integrate the S-RONS equation~\eqref{eq:srons} with the regularization parameter $\alpha = 10^{-8}$ using Matlab's \texttt{ode113s} \cite{shampine1997} with a relative and absolute error tolerance of $10^{-8}$. Figure~\ref{fig:harmonic_error_chap3} shows the relative error of the mean and covariance when using increasing number of modes $r$ in the approximate solution.
Additionally, we report computational times and errors for each simulation in Table \ref{table:harmonic_2gauss_chap3}.
We note that the cases with $r=1$ and $r=2$ are not stiff and therefore do not require any regularization, but using the same regularization parameter for every run allows for a fair comparison between all simulations. The mean is accurately predicted by RONS, regardless of how many modes we use in the approximate solution.
In this example, the true solution is a Gaussian~\cite{bruna22}, and therefore we expect to capture the true mean within the accuracy of our time integration error tolerances with our approximate solution which is a sum of Gaussians. As a result, the mean is predicted quite accurately even when $r=1$, and we do not observe a significant improvement in prediction of the mean as we use more modes. 

On the other hand, the behavior of the covariance is more complex. The true solution has a diagonal covariance matrix at the initial time. But it develops nonzero entries in the off-diagonal elements after the initial time. The approximate solution $\phat$ is the sum of Gaussians with diagonal covariance matrices, and therefore using more modes leads to significant improvement in the approximation of the covariance matrix. This is demonstrated in figure~\ref{fig:harmonic_error_chap3}, showing that as the number of modes increases from $r=1$ to $r=5$ the covariance error decreases monotonically, converging to the tolerance of the numerical time integration.
However, increasing the number of modes beyond $r=5$ does not lead to a significant improvement, indicating that 5 modes are adequate to capture the dynamics.

This example also demonstrates the scalability advantages of symbolic RONS as discussed in Section~\ref{sec:SRONS}. Namely, the number of terms $r$ can be increased without incurring additional  symbolic computational cost. That is why in Table~\ref{table:harmonic_2gauss_chap3} the symbolic computational cost is zero for $r\geq 2$. As a result, one can easily increase the number of modes until a satisfactory approximation is achieved. 

In the above examples, we assumed that the initial condition $p_0$ lies on the statistical manifold $\mathcal M$ and therefore it can be exactly represented by expansion~\eqref{eq:gen_approx} at time $t=0$. If the initial condition does not lie on the manifold initially, one can solve the optimization problem,
\begin{equation}
\btheta_0 = \argmin_{\btheta\in\Omega}\|\phat(\cdot,\btheta)-p_0\|_H^2,
\end{equation}
to find the closest point on the manifold $\mathcal M$ to the initial condition $p_0$. This needs to be solved only once at the initial time. 

%

\begin{table}
	\centering
	\begin{tabular}{|l||cccc|}
		\hline
		\textbf{Harmonic Trap} &
		\begin{tabular}[c]{@{}c@{}}Symbolic \\ computation\end{tabular} &
		\begin{tabular}[c]{@{}c@{}}Time \\ integration\end{tabular} &
		\begin{tabular}[c]{@{}c@{}}Relative error\\ of mean\end{tabular} &
		\begin{tabular}[c]{@{}c@{}}Relative error \\ of covariance\end{tabular} \\ \hline
		\multicolumn{1}{|l||}{$r=1$} &
		13.7 minutes &
		0.06 seconds &
		$\approx 5  \times 10^{-8}$ &
		$\approx 8  \times 10^{-2}$ \\ 
		\hline
		\multicolumn{1}{|l||}{$r=2$} &
		0 minutes &
		0.34 seconds &
		$\approx 10^{-8} - 10^{-9}$ &
		$\approx 2 \times 10^{-2}$\\
		\hline
		\multicolumn{1}{|l||}{$r=3$} &
		0 minutes &
		1.32 seconds &
		$\approx 10^{-8} - 10^{-9}$ &
		$\approx 2 \times 10^{-3}$\\
		\hline
		\multicolumn{1}{|l||}{$r=4$} &
		0 minutes &
		16.44 seconds &
		$\approx 10^{-8} - 10^{-9}$ &
		$\approx 2 \times 10^{-4}$\\
		\hline
		\multicolumn{1}{|l||}{$r=5$} &
		0 minutes &
		85.00 seconds &
		$\approx 10^{-8} - 10^{-9}$ &
		$\approx 10^{-5}$\\
		\hline
		\multicolumn{1}{|l||}{$r=6$} &
		0 minutes &
		470.50 seconds &
		$\approx 10^{-8} - 10^{-9}$ &
		$\approx 10^{-5}$\\
		\hline
	\end{tabular}
	\caption{Comparison of computational time and accuracy for harmonic trap example using S-RONS for various numbers of modes and regularization parameter of $\alpha = 10^{-8}$.}
	\label{table:harmonic_2gauss_chap3}
\end{table}

\section{Conclusions}\label{sec:conc_FP}
We showed that the method of reduced-order nonlinear solutions (RONS) leads to a fast and scalable method for approximating the solutions of the Fokker--Planck equation. In particular, we considered the approximate solution as the sum of shape-morphing modes, where the modes are Gaussians with time-dependent amplitudes, means and covariances. RONS equations provide a system of ODEs for optimally evolving the shape parameters such that the approximate solution stays close to a true solution of the PDE. The feasibility of approximating the probability density with a sum of Gaussians is guaranteed by the universal approximation theorem of Park and Sandberg~\cite{Park1991}.

We demonstrated the efficacy of RONS on several examples. First, we considered the Ornstein--Uhlenbeck process, where the exact solution to the corresponding Fokker--Planck equation is known. In this case, RONS reproduces this exact solution. We also considered three more complex examples, showing that RONS returns accurate approximations of the transient dynamics as well as the equilibrium density. At the same time, RONS is considerably faster than conventional methods. For instance, in the case of the stochastic Duffing oscillator, RONS is 600 times faster than direct Monte Carlo simulations.

We considered two computational methods for forming the RONS equations: symbolic RONS (or S-RONS) and collocation RONS (or C-RONS). In symbolic RONS, we use symbolic computing to
evaluate the inner products on the underlying Hilbert space $H$. This requires only 9 symbolic computation which is independent of the dimension of the system and the number of terms in the approximate solution. If the underlying Hilbert space is chosen to be the Lebesgue space $L^2$, existing symbolic computing packages easily return closed-form symbolic expressions for the required inner products. We also showed that, if we use the weighted Lebesgue space $H=L^2_\mu$ where $\mu=\phat^{-1}\id \vc x$, the metric tensor in RONS coincides with the Fisher information matrix defined on statistical manifolds. 

Our numerical experiments show that this choice of the underlying Hilbert space ($H=L^2_\mu$) leads to more accurate approximation of the true solution.
However, existing symbolic computing packages did not return a closed-form expression for the required inner products on the space $L^2_\mu$. 
In such cases, where obtaining symbolic expressions is not feasible, we used C-RONS which minimizes the error at prescribed collocation points. Consequently, C-RONS does not require computing any functional inner products and therefore is applicable in any function space. However, since it relies on collocation points, C-RONS is not scalable to higher dimensions.

Given its higher accuracy, scalability, and connection to Fisher information, using symbolic computing in the Hilbert space $L^2_\mu$ is highly desirable. Future work will explore possible avenues for incorporating symbolic computing with this weighted Lebesgue space.

\subsection*{Funding}
This work was supported by the National Science Foundation through the award DMS-2208541.


\end{document}